%% file: Lang.tex
\documentclass{article}

\usepackage{spconf}
\usepackage{graphicx}
\usepackage{amssymb}
\usepackage[cmex10]{amsmath}
\usepackage{pgfplots, filecontents}
\usepgfplotslibrary{patchplots}
\usepackage[]{algorithm2e}
\usepackage{enumerate}
\usepackage{color}
\usepackage{epstopdf}
\graphicspath{{./fig/}}
\usepackage{cite}
\newcommand{\ve}{\mathbf}
\newcommand{\m}{\mathbf}

\title{Parameter Estimation Under Model Uncertainties by Iterative Covariance Approximation}

\name{O. Lang, M. Lunglmayr and M. Huemer}
\address{Institute of Signal Processing\\
Johannes Kepler University Linz\\
Altenbergerstra{\ss}e 69, 4040 Linz, Austria}

\begin{document}
%
\maketitle
\begin{abstract}
We propose a novel iterative algorithm for estimating a deterministic but unknown parameter vector in the presence of model uncertainties. This iterative algorithm is based on a system model where an overall noise term describes both, the measurement noise and the noise resulting from the model uncertainties. This overall noise term is a function of the true parameter vector, allowing for an iterative algorithm. The proposed algorithm can be applied on structured as well as unstructured models and it outperforms prior art algorithms for a broad range of applications. 
\end{abstract}
%

\section{Introduction}
\label{sec:intro}
\input{Introduction}

\section{System Model}
\label{sec:System Model}

\input{System_Model}

\section{Iterative Algorithm}
\label{sec:Iterative Algorithm}
\input{Iterative_Algorithm}

\section{Simulation Results}
\label{sec:Simulation Results}

\input{Simulation_Results}
\section{Conclusions}
\label{sec:concl}

\input{Conclusions}

%



\bibliography{References}

\end{document}

%% file: Introduction.tex
The linear model
\begin{equation}
\ve{y}  = \m{H} \ve{x} + \ve{n} \label{equ:DECONV071}
\end{equation}
is frequently used in many areas of signal processing. Here, $\ve{y} \in \mathbb{R}^{N_\ve{y} \times 1}$ is the vector of measurements, $\ve{x} \in \mathbb{R}^{N_\ve{x} \times 1}$ is a deterministic but unknown parameter vector, $\m{H} \in \mathbb{R}^{N_\ve{y} \times N_\ve{x}}$ is the measurement matrix with $N_\ve{y} > N_\ve{x}$ and full rank, and $\ve{n} \in \mathbb{R}^{N_\ve{y} \times 1}$ is zero mean measurement noise with known covariance matrix $\m{C}_{\ve{n}\ve{n}}$. The probability density function (PDF) of $\ve{n}$ is otherwise arbitrary. Linear classical estimators such as the least squares (LS) estimator or the best linear unbiased estimator (BLUE) \cite{Kay-Est.Theory, Best_linear_unbiased_estimator_algorithm_for_received_signal-strength_based_localization} assume that the measurement matrix $\m{H}$ is perfectly known. In practice, this assumption often does not hold. A prominent case is where $\m{H}$ is a convolution matrix that is itself estimated from an imperfectly measured system output. The error in $\m{H}$ is often neglected since it is unknown.

There exist several ways to account for the errors in $\m{H}$. Two prominent algorithms that are related to the approach in this work can be found in \cite{Hierarchical_Bayesian_image_restoration}. These algorithms were derived for the task of image restoration, where the point-spread function that distorts the image is considered to be the sum of a known mean and an unknown zero-mean random part. It also provides an algorithm in the Bayesian context. In this work, however, classical estimation is considered. Hence, no prior distribution about $\ve{x}$ is assumed.


In contrast to the LS estimator and the BLUE, total least squares (TLS) estimation techniques incorporate model errors. E.g., for independent and identically distributed (i.i.d.) model errors with Gaussian PDF, the maximum likelihood (ML) solution of the TLS problem was analyzed in \cite{Linear_Regression_With_Gaussian_Model_Uncertainty_Algorithms_and_Bounds}. However, in many practical applications $\m{H}$ has some sort of structure as it is the case for Toeplitz or Hankel matrices. Then, the model errors are clearly not i.i.d. any more. Structured total least squares (STLS) techniques have been developed to deal with these kind of problems \cite{Structured_Least_Squares_Problems_and_Robust_Estimators, Formulation_and_solution_of_structured_total_least_norm_problems_for_parameter_estimation, Consistency_of_the_structured_total_least_squares_estimator_in_a}. An overview of different TLS and STLS methods can be found in \cite{Structured_total_least_squares_and_L2_approximation_problems, An_Analysis_of_the_Total_Least_Squares_Problem, Overview_of_total_least-squares_methods}. 

In this work we compare our novel approach with two iterative algorithms, which serve as performance reference in the remainder of this paper. The first one, introduced in \cite{Linear_Regression_With_Gaussian_Model_Uncertainty_Algorithms_and_Bounds}, is an approach for solving the maximum likelihood (ML) problem based on classical expectation-maximization (EM) \cite{Maximum_Likelihood_from_Incomplete_Data_via_the_EM_Algorithm}. This algorithm, referred to as ML-EM algorithm, treats the model errors as random and allows for an incorporation of the model error variance. By doing so, a uniform variance for every element in $\m{H}$ was assumed. 
The second one represents an algorithm from the class of STLS approaches and is introduced in \cite{Total_Least_Norm_Formulation_and_Solution_for_Structured_Problems}. This iterative algorithm is called the structured total least norm (STLN) algorithm and it is capable of dealing with structured measurement matrices. This approach treats the model errors as deterministic but unknown. Hence, it prevents the usage of model error variances. 

In this paper, we propose a novel iterative algorithm that incorporate information about the model error variances. Moreover, this algorithm can be employed on structured as well as unstructured problems. In contrast to the ML-EM algorithm, the algorithm is capable of incorporating different variances for every element of $\m{H}$. A difference to the STLN algorithm is that the proposed algorithm treats the model errors as random variables, allowing to incorporate the model error variances. All three algorithms require solving an inverse linear problem at each iteration.
Simulation examples are presented which show that the proposed algorithm is able to outperform both competing algorithms in a mean square error (MSE) sense for a broad range of model error and noise variances. 

The proposed iterative algorithm is based on a system model where an overall noise term describes both, the measurement noise and the noise resulting from the model uncertainties. The covariance matrix of this overall noise term is evaluated for different cases. Considering the model errors as random with known second order statistics (but otherwise arbitrary PDF) is motivated by practical examples such as multiple-input multiple-output (MIMO) communication channels or beamforming  \cite{Minimax_estimation_of_deterministic_parameters_in_linear_models_with_a_random_model_matrix, From_theory_to_practice_an_overview_of_MIMO_space-time_coded_wireless_systems, Robust_adaptive_beamforming_for_general-rank_signal_models, Non-Systematic_Complex_Number_RS_Coded_OFDM_by_Unique_Word_Prefix}.

The remainder of this paper is organized as follows: In Sec.~\ref{sec:System Model}, the underlying system model is introduced. Here we distinguish between unstructured and structured measurement matrices. For the structured case, we considered convolution matrices in this work. However, extensions to other kind of structured matrices are easily possible. The proposed iterative algorithm is discussed in Sec.~\ref{sec:Iterative Algorithm}. Simulation results demonstrating its performance are given in Sec.~\ref{sec:Simulation Results}.
\\ 
\emph{Notation:}
\\ 
Lower-case bold face variables ($\ve{a}$, $\ve{b}$,...) indicate vectors, and upper-case bold face variables ($\m{A}$, $\m{B}$,...) indicate matrices. We further use $\mathbb{R}$ and $\mathbb{C}$ to denote the set of real and complex numbers, respectively, $(\cdot)^T$ to denote transposition, $\m{I}^{n\times n}$ to denote the identity matrix of size $n\times n$, and $\m{0}^{m\times n}$ to denote the all-zero matrix of size $m\times n$. If the dimensions are clear from the context we simply write $\m{I}$ and $\m{0}$, respectively.
$E[\cdot]$ denotes the expectation operator, $[\cdot]_{i}$ the $i^\text{th}$ element of a vector and  $[\cdot]_{i,j}$ the element of a matrix at the $i^\text{th}$ row and the $j^\text{th}$ column.

%% file: System_Model.tex
This section describes the underlying model used in the remainder of this paper. In a first step, the measurement matrix is assumed to be unstructured and the model uncertainties are assumed to be independent. Afterwards, $\m{H}$ is assumed to be a structured convolution matrix built from an estimated or measured impulse response. Hence, $\m{H}$ is a special form of a Toeplitz matrix and, as it will be shown, results in correlated model uncertainties.

\subsection{Unstructured Measurement Matrices} \label{General Measurement Matrices}

We denote $\hat{\m{H}}$ as the measured or estimated measurement matrix and assume it comes along with error variances for every entry. The error variances assembled in a matrix of the same size as $\hat{\m{H}}$ is denoted as $\m{V} \in \mathbb{R}^{ N_\ve{y} \times N_\ve{x}}$. Furthermore, the errors are assumed to be independent zero mean random variables. The measurements are modeled as
\begin{equation}
\ve{y}  = \m{H}\ve{x} + \ve{n} = (\hat{\m{H}} + \m{B} ) \ve{x} + \ve{n}, \label{equ:DECONV072}
\end{equation}
where $\m{H} = \hat{\m{H}} + \m{B}$, with $\hat{\m{H}}$ being the estimated measurement matrix and $\m{B}$ being a zero mean random matrix. In \eqref{equ:DECONV072}, $\m{H}$ and $\m{B}$ are unknown while $\hat{\m{H}}$ is known. We further rewrite \eqref{equ:DECONV072} according to
\begin{eqnarray}
\ve{y}  &=&  \hat{\m{H}}\ve{x}  + \underbrace{\m{B} \ve{x} + \ve{n}}_{\ve{w}} \label{equ:DECONV044} \\
&=&  \hat{\m{H}}\ve{x}  + \ve{w}, \label{equ:DECONV045}
\end{eqnarray}
with the new overall noise vector $\ve{w}$. This noise vector combines the measurement noise with the noise from the model uncertainties. Let $\ve{b}_i^T$ be the $i^\text{th}$ row of $\m{B}$, then the $i^\text{th}$ element of $\ve{w}$ is given by
\begin{align}
[\ve{w}]_i &= \ve{b}_i^T \ve{x} + [\ve{n}]_i \label{equ:DECONV046a}
\end{align}
Since $[\ve{w}]_i$ is evaluated as the scalar product of a vector with zero mean random elements with an unknown but deterministic vector plus $[\ve{n}]_i$, $[\ve{w}]_i$ has zero mean and its variance in dependence of the unknown parameter vector $\ve{x}$ can be derived as 
\begin{align}
\sigma_i^2 =&  [\m{V}]_{i,1} |[\ve{x}]_1|^2 + [\m{V}]_{i,2} |[\ve{x}]_2|^2 + \cdots + [\m{V}]_{i,{N_\ve{x}}} |[\ve{x}]_{N_\ve{x}}|^2
 \nonumber \\
 & + [\m{C}_{\ve{n}\ve{n}}]_{i,i}.\label{equ:DECONV046b}
\end{align}
All variances assembled in a covariance matrix are combined in
\begin{equation}
\m{C}_{\ve{w}\ve{w}} = \text{diag}( \m{V} |\ve{x}|^2) + \m{C}_{\ve{n}\ve{n}}, \label{equ:DECONV048}
\end{equation}
where the term $|\ve{x}|^2$ represents a column vector of the element-wise absolute squares of the vector $\ve{x}$.

\subsection{Convolution Matrices} \label{Convolution Matrices}

We will now assume that $\m{H}$ is a linear convolution matrix constructed from the impulse response $\ve{h} \in \mathbb{R}^{ N_\ve{h} \times 1}$ of a linear system such that $\m{H}\ve{x}$ describes the convolution of the underlying sequences $h[n]$ and $x[n]$. An extension to other structured measurement matrices is easily possible. Let $\m{H}= \hat{\m{H}} + \m{B}$ have the dimension $N_\ve{y} \times N_\ve{x}$ where $N_\ve{y} = N_\ve{x} + N_\ve{h} - 1$. The $i^\text{th}$ column of the convolution matrix is defined as
\begin{equation}
[\m{H}]_{:,i} = \begin{bmatrix}
\ve{0}^{(i-1)  \times 1} \\ \ve{h} \\ \ve{0}^{(N_\ve{x} -i) \times 1}
\end{bmatrix}, \hspace{5mm}
[\hat{\m{H}}]_{:,i} = \begin{bmatrix}
\ve{0}^{(i-1)  \times 1} \\ \hat{\ve{h}} \\ \ve{0}^{(N_\ve{x} -i) \times 1}
\end{bmatrix}, \nonumber
\end{equation}
\begin{equation}
[\m{B}]_{:,i} = \begin{bmatrix}
\ve{0}^{(i-1)  \times 1} \\ \ve{e} \\ \ve{0}^{(N_\ve{x} -i) \times 1}
\end{bmatrix} \hspace{16mm}\forall i = 1,\ldots, N_\ve{x} \label{equ:DECONV057a}
\end{equation}
where $\hat{\ve{h}}$ is the estimated impulse response and $\ve{e}$ is the unknown error of $\hat{\ve{h}}$ with known error covariance matrix $\m{C}_{\ve{e}\ve{e}}\in \mathbb{R}^{N_\ve{h} \times N_\ve{h}}$. In this case, the model uncertainties of $\hat{\m{H}}$ are clearly not independent anymore, leading to a different calculation of $\m{C}_{\ve{w}\ve{w}}$.

 Let $\ve{b}'_i = [\m{B}]_{:,i}$ denote the $i^\text{th}$ column of $\m{B}$. The subsequent column $\ve{b}'_{i+1}$ can be derived by shifting down the elements of $\ve{b}'_i$ by one position:
 \begin{equation}
 \ve{b}'_{i+1} = \begin{bmatrix}
 \ve{0}^{1 \times (N_\ve{y} - 1)} & 0 \\
 \m{I}^{(N_\ve{y}-1) \times (N_\ve{y}-1)} & \ve{0}^{(N_\ve{y}-1) \times 1}
 \end{bmatrix} \ve{b}'_i = \m{D}\ve{b}'_i.  \label{equ:DECONV057}
 \end{equation} 
 With that, the product  $\m{B} \ve{x}$ in \eqref{equ:DECONV044} follows to
  \begin{align}
  \m{B} \ve{x} =& \ve{b}'_1 [\ve{x}]_1 + \ve{b}'_2 [\ve{x}]_2 + \hdots + \ve{b}'_{N_\ve{x}} [\ve{x}]_{N_\ve{x}} \label{equ:DECONV058} \\
  		=& \left([\ve{x}]_1 \m{I} + [\ve{x}]_2 \m{D} + \hdots + [\ve{x}]_{N_\ve{x}} \m{D}^{{N_\ve{x}}-1}\right)\ve{b}'_1 \label{equ:DECONV061} \\
  		=& \m{P}(\ve{x})\ve{b}'_1. \label{equ:DECONV062}
  \end{align}
   With this result, $\ve{w}$ can be written as $\ve{w}= \m{P}(\ve{x})\ve{b}'_1 + \ve{n}$
   with the covariance matrix
   \begin{align}
   \m{C}_{\ve{w}\ve{w}} =& E\left[\left(\m{P}(\ve{x})\ve{b}'_1 \right)\left(\m{P}(\ve{x})\ve{b}'_1 \right)^H    \right] + \m{C}_{\ve{n}\ve{n}} \label{equ:DECONV065} \\
   	=& \m{P}(\ve{x}) \m{C}_{\ve{b}'_1\ve{b}'_1} \m{P}(\ve{x})^H + \m{C}_{\ve{n}\ve{n}}.\label{equ:DECONV066} 
   \end{align}
   The covariance matrix $\m{C}_{\ve{b}'_1\ve{b}'_1}$ follows from \eqref{equ:DECONV057a} and the covariance matrix of the estimation error $\ve{e}$ according to
   \begin{equation}
   \m{C}_{\ve{b}'_1\ve{b}'_1} = \begin{bmatrix}
   \m{C}_{\ve{e}\ve{e}} & \m{0}^{N_\ve{h} \times (N_\ve{x}-1)} \\
   \m{0}^{ (N_\ve{x}-1) \times N_\ve{h} } & \m{0}^{ (N_\ve{x}-1) \times (N_\ve{x}-1)}
   \end{bmatrix}\in \mathbb{R}^{N_\ve{y} \times N_\ve{y}}.
   \end{equation}
Note that this formulation allows for two sources of correlations. The first source comes from the structure in $\m{H}$. The second source of correlation comes from $\m{C}_{\ve{e}\ve{e}}$, which describes the errors in $\hat{\ve{h}}$. Hence, the iterative algorithm introduced in the next section is capable of dealing with both kind of correlations.

%% file: Iterative_Algorithm.tex
 An ideal but theoretical estimator is the BLUE applied on the linear model in \eqref{equ:DECONV072} using the \emph{true} $\m{H}$ according to
\begin{equation}
 \hat{\ve{x}} = \left( \m{H}^H  \m{C}_{\ve{n}\ve{n}}^{-1} \m{H} \right)^{-1} \m{H}^H  \m{C}_{\ve{n}\ve{n}}^{-1} \ve{y}.    \label{equ:DECONV048atwoo}
\end{equation}
This theoretical estimator is referred to as BLUE \emph{with perfect model knowledge}. Similarly, the BLUE applied on the linear model in \eqref{equ:DECONV045}, incorporating the estimated measurement matrix $\hat{\m{H}}$ but the \emph{true} covariance matrix $\m{C}_{\ve{w}\ve{w}}$ follows as
\begin{equation}
 \hat{\ve{x}} = ( \hat{\m{H}}^H  \m{C}_{\ve{w}\ve{w}}^{-1} \hat{\m{H}} )^{-1} \hat{\m{H}}^H  \m{C}_{\ve{w}\ve{w}}^{-1} \ve{y}    \label{equ:DECONV048a}
\end{equation}
and is referred to as BLUE \emph{with perfect knowledge of} $\m{C}_{\ve{w}\ve{w}}$ \cite{Approximate_best_linear_unbiased_channel_estimation}.  The determination of the true $\m{C}_{\ve{w}\ve{w}}$ according to \eqref{equ:DECONV048} or \eqref{equ:DECONV066}, however, requires the knowledge of the true parameter vector. To overcome this problem, we propose the iterative algorithm described below. Its basic idea is to make an initial guess of the parameter vector termed $\hat{\ve{x}}_0$ (the index denotes the algorithm's iteration number). This first guess could, e.g., origin from an LS estimation which does not incorporate any noise statistics. $\hat{\ve{x}}_0$ is then used to estimate $\hat{\m{C}}_{\ve{w}\ve{w},0}$ in \eqref{equ:DECONV048} or \eqref{equ:DECONV066}. This estimated covariance matrix is then incorporated by the BLUE in order to yield a better estimate $\hat{\ve{x}}_1$ and so on. This procedure is summarized as shown in Algorithm \ref{alg:Proposed_algorithm}.

\begin{algorithm}[h]
\SetKwInOut{Initialization}{Initialization}
 \Initialization{LS estimation}
 $\hspace{5mm} \hat{\ve{x}}_0 = \left(\hat{\m{H}}^T \hat{\m{H}} \right)^{-1} \hat{\m{H}}^T  \ve{y}$\;
 \For{$k\leftarrow 0$ \KwTo $N_\text{iter}$}{
 estimate $\m{C}_{\ve{w}\ve{w},k}$ according to \eqref{equ:DECONV048} or \eqref{equ:DECONV066}  using $\hat{\ve{x}}_k$ instead of $\ve{x}$ \;
  $\hat{\ve{x}}_{k+1} =  \left(\hat{\m{H}}^T \hat{\m{C}}_{\ve{w}\ve{w},k}^{-1}\hat{\m{H}} \right)^{-1} \hat{\m{H}}^T \hat{\m{C}}_{\ve{w}\ve{w},k}^{-1} \ve{y}$ \; 
 } 
  \caption{proposed algorithm}
  \label{alg:Proposed_algorithm}
\end{algorithm}

The proposed algorithm is of similar complexity as the ML-EM and STLN algorithms. It performs a weighting of the measurements according to $\hat{\m{C}}_{\ve{w}\ve{w},k}$, which incorporates the model error variances as well as the measurement noise variances. In the case of $\m{H}$ being a convolution matrix, even the covariances of the estimated impulse response are considered in order to improve the estimation. 

Note that for both cases $\hat{\m{C}}_{\ve{w}\ve{w},k}$ is almost surely invertible since $\m{C}_{\ve{n}\ve{n}}$ serves as a regularization term in \eqref{equ:DECONV048} and \eqref{equ:DECONV066}.


Although convergence cannot be ensured, simulations showed that divergence is a rare exception for reasonable values of $\m{V}$.

A stopping criteria can be implemented in several ways. One possibility is to stop the iterations when $\hat{\ve{x}}$ does not significantly change from one iteration to the next. Simulations showed that the major performance gain is usually achieved after the first iteration. Hence, a predefined number of iterations may be utilized instead of a stopping criteria. 

Naturally, there exists at least one case where the iterations yield no performance gain. If $\hat{\m{C}}_{\ve{w}\ve{w},k}$ is a scaled identity matrix, the proposed algorithm reduces to the ordinary LS estimator, preventing any performance increase. This is, e.g., the case when the following two conditions hold: a) The measurement matrix is unstructured and $\m{V}$ has the same variance at every element. b) the noise covariance matrix $\m{C}_{\ve{n}\ve{n}}$ is a scaled identity matrix. 

We note that a similar iterative application of the BLUE was applied in \cite{Pilot_Assisted_Time_Varying_Channel_Estimation_for_OFDM_Systems, Approximate_best_linear_unbiased_channel_estimation, Improved_channel_estimation_using_superimposed_training} for channel impulse response estimation in wireless communication applications. Compared to them, the proposed algorithm is applicable to various applications with structured or unstructured model uncertainties. In \cite{The_bias_of_the_unbiased_estimator_A_study_of_the_iterative_application_of_the_BLUE_method} investigations of a similar procedure as the presented algorithm can be found but only for a very simplified model compared to the investigations in this work. As a result of that, the algorithms presented in \cite{Pilot_Assisted_Time_Varying_Channel_Estimation_for_OFDM_Systems, Approximate_best_linear_unbiased_channel_estimation, Improved_channel_estimation_using_superimposed_training, The_bias_of_the_unbiased_estimator_A_study_of_the_iterative_application_of_the_BLUE_method} are not considered in the following simulations. Here, we rather compare the proposed algorithm with the STLS algorithm \cite{Total_Least_Norm_Formulation_and_Solution_for_Structured_Problems}, the ML-EM algorithm \cite{Linear_Regression_With_Gaussian_Model_Uncertainty_Algorithms_and_Bounds} as well as the estimators in \eqref{equ:DECONV048atwoo} and \eqref{equ:DECONV048a}.

%% file: Simulation_Results.tex
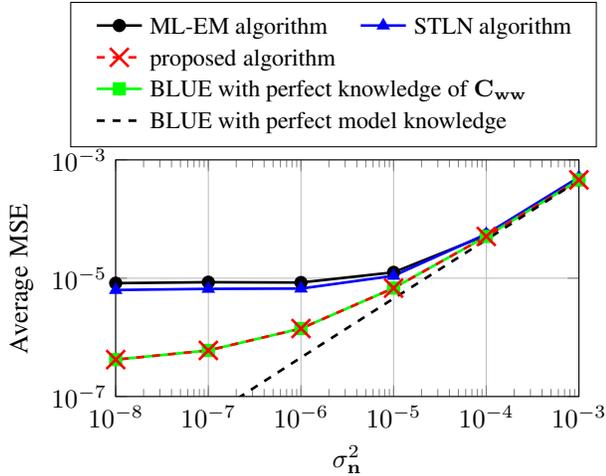
\begin{figure}[t]
\begin{center}
\begin{tikzpicture}
\begin{loglogaxis}[compat=newest, 
width=0.9\columnwidth, height = .55\columnwidth, xlabel=$\sigma_\ve{n}^2$, 
ylabel style={align=center}, 
ylabel style={text width=3.4cm},
ylabel={Average MSE}, 
legend pos=north east, 
legend cell align=left,
legend columns=2, 
        legend style={
            /tikz/column 2/.style={
                column sep=5pt,
            },
        font=\small},
xmin = 0.00000000999,
xmax = 0.001,
ymax = 0.001,
ymin = 0.999e-7,
grid=major,
legend style={
at={(-0.05,1.5)},
anchor=north west}
]

\addplot[line width=1pt, color=black, mark=otimes*, mark options={solid},] table[x index =0, y index =6] {ICASSP_Convmtx_1.dat};
\label{p1}

\addplot[line width=1pt, color=blue, style=solid, mark=triangle*, mark options={solid}] table[x index =0, y index =4] {ICASSP_Convmtx_1.dat};
\label{p2}

\addplot[line width=1pt, color=green, mark=square*, mark options={solid}] table[x index =0, y index =3] {ICASSP_Convmtx_1.dat};
\label{p3}

\addplot[line width=1pt, color=red, style=dashed, mark=x, mark options={solid}, mark size=5pt] table[x index =0, y index =2] {ICASSP_Convmtx_1.dat};
\label{p4}

\addplot[line width=1pt, color=black, style=dashed] table[x index =0, y index =7] {ICASSP_Convmtx_1.dat};
\label{p5}

\end{loglogaxis}

\node [draw,fill=white,anchor=north east] at (6.4,5.25) {\shortstack[l]{
	\begin{tabular}{ll}
		\ref{p1} \small{ML-EM algorithm} &  \ref{p2} \small{STLN algorithm} \\
		\multicolumn{2}{l}{\ref{p4} \small{proposed algorithm}} \\
		\multicolumn{2}{l}{\ref{p3} \small{BLUE with perfect knowledge of $\m{C}_{\ve{w}\ve{w}}$}} \\  
		\multicolumn{2}{l}{\ref{p5} \small{BLUE with perfect model knowledge}} 
  	\end{tabular} }};

\end{tikzpicture}
\caption{ Average MSEs of different iterative algorithms plotted over the noise variance $\sigma_\ve{n}^2$.  \label{fig:Example1_Fig1} }
\end{center}
\end{figure}

\begin{figure}[t]
\begin{center}
\begin{tikzpicture}
\begin{semilogyaxis}[compat=newest, 
width=0.9\columnwidth, height = .55\columnwidth, xlabel=Iteration, 
ylabel style={align=center}, 
ylabel style={text width=3.4cm},
ylabel={Average MSE}, 
legend pos=north east, 
legend cell align=left,
legend columns=2, 
        legend style={
            /tikz/column 2/.style={
                column sep=5pt,
            },
        font=\small},
xmin = 0,
xmax = 10,
ymax = 0.00001,
ymin = 0.000001,
grid=major,
legend style={
at={(-0.05,1.37)},
anchor=north west}
]

\addplot[line width=1pt, color=black, mark=otimes*, mark options={solid},] table[x index =0, y index =3] {ICASSP_Convmtx_conv_curves.dat};
\addlegendentry{{ML-EM algorithm}}

\addplot[line width=1pt, color=blue, style=solid, mark=triangle*, mark options={solid}] table[x index =0, y index =2] {ICASSP_Convmtx_conv_curves.dat};
\addlegendentry{{STLN algorithm}}

\addplot[line width=1pt, color=red, style=dashed, mark=x, mark options={solid}, mark size=5pt] table[x index =0, y index =1] {ICASSP_Convmtx_conv_curves.dat};
\addlegendentry{{proposed algorithm}}

\end{semilogyaxis}
\end{tikzpicture}
\caption{Average MSE values plotted over the number of iterations. 
 \label{fig:Example1_conv_plot} }
\end{center}
\end{figure}
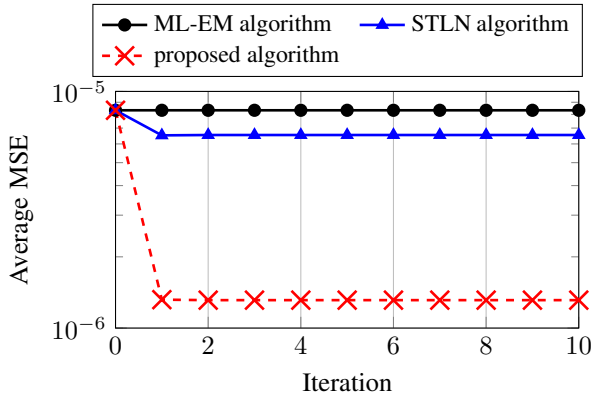

In this example, $\m{H} \in \mathbb{R}^{7 \times 3}$ is a convolution matrix and describes the discrete convolution of the impulse response $h[n]$ with signal $x[n]$. The vector notations of $h[n]$ and $x[n]$ are given by $\ve{h} \in \mathbb{R}^{5 \times 1}$ and $\ve{x} \in \mathbb{R}^{3 \times 1}$, respectively. For the simulations, the impulse responses is randomly generated from a Gaussian distribution with mean $E[\ve{h}] = \ve{0}^{5 \times 1}$ and covariance matrix $\m{C}_{\ve{h}\ve{h}}= \m{I}^{5 \times 5}$. The input signal is chosen to be $\ve{x} = \begin{bmatrix}
1 & 0.5 & 0.25 \end{bmatrix}^T$. For the first analysis, the noise covariance matrix is a scaled identity matrix $\m{C}_{\ve{n}\ve{n}}= \sigma_\ve{n}^2 \m{I}^{7 \times 7}$, where the scaling factor $\sigma_\ve{n}^2$ is varied between $10^{-8}$ and $10^{-3}$. The impulse response estimation step is assumed to yield zero mean errors with error covariance matrix $\m{C}_{\ve{e}\ve{e}}= \text{diag}\left(\begin{bmatrix} 10^{-4} & 10^{-5} & 10^{-6} & 10^{-6} & 10^{-6} \end{bmatrix}\right)$. 
For this model, the proposed algorithm in Sec.~\ref{sec:Iterative Algorithm} is compared with the BLUE with perfect model knowledge in \eqref{equ:DECONV048atwoo}, the BLUE with perfect knowledge of $\m{C}_{\ve{w}\ve{w}}$ in \eqref{equ:DECONV048a}, the ML-EM algorithm, and the STLN algorithm. For the latter one the \emph{$l_2$} norm minimization, a tolerance $\epsilon = 10^{-10}$ and $\m{D} = \m{I}^{5 \times 5}$ is chosen. Furthermore, $\m{X}$ ( (2.1) in \cite{Total_Least_Norm_Formulation_and_Solution_for_Structured_Problems}) is identified to be the first $N_\ve{h}$ columns of $\m{P}(\ve{x})$ in \eqref{equ:DECONV062}. For more details on these parameters we refer to \cite{Total_Least_Norm_Formulation_and_Solution_for_Structured_Problems}. For the ML-EM algorithm $\sigma_\ve{h}^2$ is set to the mean value of $\m{V}$ \cite{Linear_Regression_With_Gaussian_Model_Uncertainty_Algorithms_and_Bounds}. While the STLN algorithm comes with its own termination criterium for which we choose $\epsilon = 10^{-10}$ \cite{Total_Least_Norm_Formulation_and_Solution_for_Structured_Problems}, the proposed algorithm and the ML-EM algorithm were executed for $N_\text{iter}=10$ iterations. However, as we discuss below, $N_\text{iter}$ could be reduced significantly. The resulting MSE values averaged over the elements of the MSE vector are presented in Fig.~\ref{fig:Example1_Fig1}. This figure shows that the proposed algorithm attains the performance bound given by the BLUE with perfect knowledge of $\m{C}_{\ve{w}\ve{w}}$ and outperforms the competing algorithms especially for low $\sigma_\ve{n}^2$. The performance gain is more than one order of magnitude in MSE for small noise variances. For large noise variances all investigated algorithms perform approximately equal. The reason for this is that the model uncertainties vanish compared to the large measurement noise samples in that case. For the same reason, the gap between all considered algorithms and the BLUE with perfect model knowledge decreases with increasing noise variance. Simulations showed that, if one would have chosen $\m{C}_{\ve{e}\ve{e}}$ to be a scaled identity matrix, the STLN algorithm would have similar performance as the proposed algorithm for very low noise variances. Furthermore, simulations showed that the performance gain approximately stays the same for other values of $\ve{x}$. Fig.~\ref{fig:Example1_conv_plot} shows the convergence behavior of the algorithms for $\sigma_\ve{n}^2 = 10^{-6}$. First of all, one recognizes that the ML-EM algorithm is not able to significantly improve the estimation accuracy compared to the initial LS estimation in this example. Furthermore, it shows that the STLN algorithm as well as the proposed algorithm achieve most of their performance gain during the first iteration. This extremely fast convergence allows to reduce the number of iterations to one without any significant loss in performance.

%% file: Conclusions.tex
In this work, a novel iterative algorithm for estimating an unknown but deterministic parameter vector in the presence of model errors and measurement noise is presented. This algorithm iteratively estimates the covariance matrix of an overall noise term, which describes the effects of the measurement noise as well as the noise resulting from the model uncertainty. This overall noise term was analyzed for unstructured model errors and for the case where the measurement matrix is a convolution matrix. For the latter case, simulation results are presented demonstrating the performance gain compared to competing algorithms. Convergence curves demonstrate the extremely fast convergence of the proposed algorithm, which achieves almost its optimum estimation accuracy after a single iteration.